\documentclass[notitlepage,leqno,11pt]{article}
\usepackage{amssymb}
\catcode`\@=11
\@addtoreset{equation}{section}

\catcode`\@=12

\usepackage{latexsym}
\usepackage{amsmath}
\usepackage{amsfonts}
\usepackage{amssymb}
\usepackage{mathrsfs}

\renewcommand{\d}{\delta }

\newcommand{\e }{\varepsilon }

\renewcommand{\l }{\lambda }
\newcommand{\s }{\sigma }

\renewcommand{\O }{\Omega }
\newcommand{\be}{\begin{equation}}
\newcommand{\ee}{\end{equation}}

\newcommand{\R}{\mathbb{R}}
\newcommand{\de}{\partial}
\def\C{{\mathcal C}}
\def\S{{\mathbb S}}
\def\f{{\varphi}}

\def\eps{{\varepsilon}}

\def\weak{{\,\rightharpoonup\,}}

\newcommand{\disc }{\mathbb{D}}

\newtheorem{Theorem}{Theorem}[section]
\newtheorem{Lemma}[Theorem]{Lemma}
\newtheorem{Proposition}[Theorem]{Proposition}
\newtheorem{Corollary}[Theorem]{Corollary}
\newtheorem{Remark}[Theorem]{Remark}

\def\proof{\noindent{{\bf Proof. }}}
\def\square{\vbox{
    \hrule height .4pt
    \hbox{\vrule width .4pt height 7pt \kern 7pt
       \vrule width .4pt}
    \hrule height .4pt }}

\def\QED{\hfill {$\square$}\goodbreak \medskip}

\begin{document}

\title{ Hardy-Poincar\'e inequalities with  boundary
singularities}

\author{Mouhamed Moustapha Fall\footnote{\footnotesize{Universit\'e Chatolique de Louvain-La-Neuve,
d\'epartement de math\'ematique.
 Chemin du Cyclotron 2,  1348 Louvain-la-Neuve, Belgique.
E-mail: {\tt mouhamed.fall@uclouvain.be.}  Work partially supported
by FIRB �"Analysis and Beyond", 2009-2012.} } ~and Roberta
Musina\footnote{\footnotesize{Dipartimento di Matematica ed
Informatica, Universit\`a di Udine, via delle Scienze, 206 -- 33100
Udine, Italy. E-mail: {\tt musina@dimi.uniud.it.}
}}}
\date{}

\maketitle

\bigskip

\noindent {\footnotesize{\bf Abstract.}  Let $\O$ be a bounded
domain in $\R^N$ with $0\in\de\O$ and $N\ge 2$. In this paper we
study the Hardy-Poincar\'e inequality for maps in $H^1_0(\Omega)$.
In particular we give sufficient and some necessary conditions
so that the best constant is achieved.}
\bigskip\bigskip



\bigskip
\bigskip
\vfill\eject


\section*{Introduction}

Let $\Omega$ be a smooth bounded domain in $\R^N$, with  $N\ge 2$. In this
paper we assume that $0\in\partial \Omega$ and we study the
minimization problem
\begin{equation}
\label{eq:problem}
\mu_{\l}(\O):=
\inf_
{\genfrac{}{}{0pt}{}{\scriptstyle{u\in H^1_0(\O)}}{\scriptstyle{u\neq 0}}}
~\frac{\displaystyle\int_{\O}|\nabla u|^2~dx-\l\int_{\O}|u|^2~dx}
{\displaystyle\int_{\O}|x|^{-2}|u|^2~dx}~,
\end{equation}
where $\lambda\in\R$ is a varying parameter. For $\lambda=0$ the
{\em $\Omega$-Hardy constant} $\mu_0(\O)\ge (N-2)^2/4$ is the best constant in
the Hardy inequality for maps supported by $\Omega$. If $N=2$ it has been proved in \cite{CaMuUMI},
Theorem 1.6,
that $\mu_0(\O)$ is positive.

Problem (\ref{eq:problem}) carries some similarities with the questions studied by Brezis and Marcus
in \cite{BM}, where the weight is
the inverse-square of  the distance from the boundary of $\Omega$. Also the paper
\cite{DD} by D\'avila and Dupaigne is somehow related to the minimization problem (\ref{eq:problem}). Indeed,
notice that for any fixed $\lambda\in\R$, any extremal for $\mu_\lambda(\O)$
is a weak solution to the linear Dirichlet problem
\begin{equation}
\label{eq:equation}
\begin{cases}
-\Delta u=\mu |x|^{-2} u+\lambda u\quad
&\textrm{on~ $\Omega  $}\\
u=0&\textrm{on~ $ \partial\Omega $,}
\end{cases}
\end{equation}
where $\mu=\mu_\lambda(\Omega)$. If $\mu_\lambda(\O)$ is achieved, then
$\mu_\lambda(\O)$ is the first eigenvalue of the operator
$-\Delta-\lambda$ on $H^1_0(\O)$. Starting from a different point of view,
for $0\in\O$, $N\ge 3$ and $\mu\le (N-2)^2/4$,
D\'avila and Dupaigne  have proved in \cite{DD} the existence of the first eigenfunction
$\f_1$ of the operator $-\Delta-\mu|x|^{-2}$ on a suitable functional space
$H(\O)\supseteq H^1_0(\O)$. Notice
that $\f_1$ solves (\ref{eq:equation}), where the eigenvalue
$\lambda$ depends on the datum $\mu$.

The problem of the existence of extremals for the $\O$-Hardy constant $\mu_0(\O)$
was already discussed in \cite{CaMuUMI} in case $N=2$ (with $\Omega$ possibly unbounded or
singular at $0\in\partial \Omega$) and in \cite{PT}, where $\O$ is a suitable compact
perturbation of a cone in $\R^N$.
Hardy-Sobolev
inequalities with singularity at the boundary have been studied by
several authors. We quote for instance  \cite{CaMuPRSE}, \cite{DN}, \cite{Eg}, \cite{GK}, \cite{GR1}, \cite{GR2}, and references there-in.

\medskip

The minimization problem (\ref{eq:problem}) is not compact, due to
the group of dilations in $\R^N$. Actually it might happen that all minimizing sequences
concentrate at $0$. In this case $\mu_\lambda(\O)$ is not achieved and $\mu_\lambda(\O)=\mu^+$, where
$$
\mu^+=\frac{N^2}{4}
$$
is the best constant in the Hardy inequality for maps with support
in a half-space. Indeed in Section \ref{S:existence}  we first show
that
\begin{equation}
\label{eq:large} \sup_{\l\in\R}\mu_\lambda(\Omega)= \mu^+~\!
\end{equation}
then we deduce that, provided $\mu_\l(\O)<\mu^+$, every minimizing
sequence for $\mu_\l(\O)$ converges in $H^1_0(\O)$ to an extremal
for $\mu_\lambda(\O)$.

We recall that $\O$ is said to be locally concave at $0\in\de\O$  if
it contains a half-ball. That is there exists  $r>0$ such that
\begin{equation}
\label{eq:halfball}
  \{x\in\R^N~|~x\cdot\nu>0\}\cap B_r(0)\subset \O~\!,
\end{equation}
where $\nu$ is the interior normal  of $\de\O$ at 0. Notice that if
all the principal curvatures of $\de\O$ at $0$, with respect to
$\nu$, are strictly negative, then condition (\ref{eq:halfball}) is
satisfied.

 Our first main
result is stated in the following theorem.
\begin{Theorem}\label{th:ext-nexJl0}
  Let $\O\in\R^N$ be a smooth bounded domain with $0\in\partial\O$.
  Assume that $\O$ is locally concave at $0$.
Then $\mu_\l(\O)$  is attained if and only if  $\mu_\l(\O)<\mu^+$.
\end{Theorem}

The "only
if" part, which is the most intriguing,  is a consequence of Corollary
\ref{C:ne} in Section \ref{S:ne}, where we provide local
nonexistence results for problem
\begin{equation}
\label{eq:inequality}
\begin{cases}
-\Delta u\ge \mu|x|^{-2}u+\lambda u\quad
&\textrm{on~\! $\Omega  $}\\
u\ge0&\textrm{in~ $\O$,}
\end{cases}
\end{equation}
also for negative values of the parameter $\lambda$.

Up to now several questions concerning the infimum $\mu_\lambda(\O)$
 are  still open. Put
\begin{equation}
\label{eq:lambda*}
\lambda^*:=\inf\left\{\lambda\in\R~|~\mu_\lambda(\O)<\mu^+~\right\}~\!.
\end{equation}
Since the map $\lambda\mapsto\mu_\lambda(\O)$ is non increasing, then
$\mu_\lambda(\O)$ is achieved for any $\lambda>\lambda^*$, by the existence
Theorem \ref{th:extJl0}. If $\l^*\in\R$, then from  (\ref{eq:large}) it follows that
$\mu_{\l}(\O)=\mu^+$ for any $\l\leq\l^*$ and hence
$\mu_\lambda(\O)$ is not achieved if $\lambda<\lambda^*$.  We don't know if there exist
domains $\O$ for which $\lambda^*=-\infty$. On the other hand
we are able to prove the following facts
(see Section \ref{S:estimates} for the precise statements):

\begin{description}
\item$i)$ If $\O$ is locally convex at $0$, that is, if there exists $r>0$ such that
$\O\cap B_r(0)$ is contained in a half-space, then $\lambda^*>-\infty$.
\item$ii)$ If $\O$ is contained in a half-space then
\begin{equation}
\label{eq:lambda*BV}
\lambda^*\geq \frac{\lambda_1(\disc)}{|\textrm{diam}(\O)|^2}~,
\end{equation}
where $\lambda(\disc)$ is the first Dirichlet eigenvalue of the unit ball
$\disc$ in $\R^2$ and $\textrm{diam}(\O)$ is the diameter of $\O$.
\item$iii)$ For any $\delta>0$ there exists $\rho_\delta>0$ such that, if
$$
\O\supseteq \{x\in\R^N~|~x\cdot\nu>-\delta|x|~,~\alpha<|x|<\beta~\}
$$
for some $\nu\in\S^{N-1}$, $\beta>\alpha>0$ with $\beta/\alpha>\rho_\delta$, then
$\lambda^*<0$. In particular the Hardy constant $\mu_0(\O)$ is achieved.
\end{description}

The relevance of the geometry of $\O$ at the origin is confirmed by
Theorem \ref{th:ext-nexJl0}, by $i)$ and by the existence theorems
proved in  \cite{GK}, \cite{GR1} and \cite{GR2} for a related
superlinear problem. However, it has to be noticed that also the
(conformal) "size" of $\O$ (even far away from the origin) has some
impact on the existence of compact minimizing sequences. Actually,
no requirement on the curvature of $\O$ at $0$ is needed in $iii)$.
In particular, there exist smooth domains having strictly positive
principal curvatures at $0$, and such that the Hardy constant
$\mu_0(\O)$ is achieved.

\bigskip

The paper is organized as follows.

In Section \ref{S:preliminaries} we point out few remarks on the
Hardy inequality  on  dilation-invariant domains.

In Section \ref{S:existence},   Theorem \ref{th:extJl0}, we give
sufficient conditions for the existence of minimizers for (\ref{eq:problem}).

In Section \ref{S:ne} we prove some  nonexistence theorems
for solutions to  (\ref{eq:inequality}) that might have an independent interest.

To prove inequality (\ref{eq:lambda*BV})  in case $\O$ is contained in a half
space, in  Section \ref{S:reminder} we provide
computable remainder terms for the  Hardy  inequality on half-balls.
We adopt here an argument by Brezis-V\'{a}zquez
\cite{BV}, where bounded domains $\O\subset \R^N$, with $N\ge 3$ and $0\in\O$,
 are considered.

In Section \ref{S:estimates} we estimate $\lambda^*$ from below and form above, under suitable assumptions on $\O$.

\bigskip
\small
\noindent
{\bf Notation}
\begin{description}
\item$\bullet$ $\R^N_+$ and $\S^{N-1}_+$ denote any half space and any hemisphere,
respectively. More precisely,
$$
\R^N_+=\{x\in\R^N~|~x\cdot\nu>0~\}~,\quad \S^{N-1}_+=\S^{N-1}\cap\R^N_+
$$
where $\nu$ is any unit vector in $\R^N$.

\item$\bullet$ $B_R(x)$ is the open ball in $\R^N$ of radius $r$ centered at $x$. If
$x=0$ we simply write $B_R$. If $N=2$ we shall often write $\disc_R$ and $\disc$
instead of $B_R$, $B_1$, respectively.

\item$\bullet$ We denote by $H^1(\S^{N-1})$ the standard
Sobolev space of maps on the unit sphere and by
$\nabla_\sigma$, $\Delta_\sigma$ the gradient and the Laplace-Beltrami operator on $\S^{N-1}$,
respectively.

\item$\bullet$  Let
$\Sigma$ be a domain in $\S^{N-1}$. We denote by $H^1_0(\Sigma)$ the closure of
$C^\infty_c(\Sigma)$ in the $H^1(\S^{N-1})$-space and by
$\lambda_1(\Sigma)$ the fist Dirichlet eigenvalue on $\Sigma$.

\item$\bullet$ For any domain $\O\subset\R^N$, we denote by
$L^2(\O;|x|^{-2}~dx)$ the space of measurable maps on $\O$ such that
$\displaystyle{\int_\O|x|^{-2}|u|^2~dx<\infty}$. We put also
$$
\widehat H^1(\O):=H^1(\O)\cap L^2(\O;|x|^{-2}~dx)~\!,
$$
where
$H^1(\O)$ is the standard Sobolev space of maps on $\O$.
\end{description}

\normalsize

\section{Preliminaries}
\label{S:preliminaries}

In this section we collect a few remarks on the Hardy inequality
on  dilation-invariant domains that are partially contained for example
in \cite{CaMuUMI} (in case $N=2$)
and in \cite{PT}.

Via polar coordinates, to any domain $\Sigma$ in $\S^{N-1}$ we
associate a {cone} $\C_\Sigma\subset\R^{N-1}$  and a (half)
{cylinder} $\mathcal Z_\Sigma\subset\R^{N+1}$ by setting
$$
\C_\Sigma:=\{~t\sigma~|~t>0~,~\sigma\in\Sigma~\}~,\quad \mathcal
Z_\Sigma:=\R_+\times\Sigma~.
$$
If $\Sigma$ is a smooth domain in $\S^{N-1}$,  then
$\C_\Sigma$ is a Lipschitz, dilation-invariant
 domain in $\R^{N-1}$.
In particular, if $\Sigma$ is a half-sphere, then $\C_\Sigma$ is a
half-space. The map
$$
\R^{N-1}\setminus\{0\}\to\R^{N+1}~,\quad x\mapsto \left(-\log|x|,\frac{x}{|x|}\right)
$$
is an homeomorphism $\C_\Sigma\to \mathcal Z_\Sigma$. It induces the
Emden-Fowler transform
$$
T:C^\infty_c(\C_\Sigma)\to C^\infty_c(\mathcal Z_\Sigma)~,\quad
u(x)=|x|^{\frac{2-N}{2}}~(Tu)\left(-\log|x|,\frac{x}{|x|}\right)~\!.
$$
A direct computation based on
the divergence theorem gives
\begin{equation}
\label{eq:EF1}
 \int_{\C_\Sigma}|\nabla u|^2~dx= \frac{(N-2)^2}{4}
 \int_0^\infty\!\!\!\!\int_\Sigma|Tu|^2~dsd\sigma+
 \int_0^\infty\!\!\!\!\int_\Sigma|\nabla_{s,\sigma}Tu|^2~dsd\sigma
\end{equation}
\begin{equation}
\label{eq:EF2}
\int_{\C_\Sigma}|x|^{-2} |u|^2~dx=  \int_0^\infty\!\!\!\!\int_\Sigma|Tu|^2~dsd\sigma,
\end{equation}
where $\nabla_{s,\sigma}=(\partial_s,\nabla_\sigma)$ denotes the gradient on $\R_+\times\S^{N-1}$.

Now we introduce the Hardy constant on the cone $\C_\Sigma$:
\begin{equation}
\label{eq:hardycone}
\mu_0(\C_\Sigma):=
\inf_
{\genfrac{}{}{0pt}{}{\scriptstyle{u\in C^\infty_c(\C_\Sigma)}}{\scriptstyle{u\neq 0}}}
~\frac{\displaystyle\int_{\C_\Sigma}|\nabla u|^2~dx}
{\displaystyle\int_{\C_\Sigma)}|x|^{-2}|u|^2~dx}~\!.
\end{equation}
In the next proposition we notice that the Hardy inequality on
$\C_\Sigma$ is equivalent to the Poincar\'e inequality for  maps
supported be the cylinder $\mathcal Z_\Sigma$.

\begin{Proposition}
\label{P:mu}
Let $\C_\Sigma$ be a cone.
Then
$$
\displaystyle{\mu_0(\C_\Sigma)}=\frac{(N-2)^2}{4}+\lambda_1(\Sigma).
$$
\end{Proposition}

\proof
By (\ref{eq:EF1}), (\ref{eq:EF2}) it turns out that
$$
\mu_0(\C_\Sigma) - \frac{(N-2)^2}{4}=\inf_
{\genfrac{}{}{0pt}{}{\scriptstyle{v\in C^\infty_c(\mathcal Z_\Sigma)}}{\scriptstyle{v\neq 0}}}
\frac{\displaystyle \int_0^\infty\!\!\!\!\int_\Sigma|\nabla_{s,\sigma} v|^2~dsd\sigma}
{\displaystyle \int_0^\infty\!\!\!\!\int_\Sigma|v|^2~dsd\sigma}=:\lambda_1(\mathcal Z_\Sigma)~\!.
$$
The result follows by noticing that $\lambda_1(\mathcal Z_\Sigma)=\lambda_1(\Sigma)$.
\QED

The eigenvalue $\lambda_1(\Sigma)$ is explicitly known in
few cases. For example, if $\Sigma=\S^{N-1}_+$ is a half-sphere then
$\lambda_1(\S^{N-1}_+)=N-1$. Thus,   the Hardy
constant of a half space is given by
\begin{equation}
\label{eq:half-space} \mu_0(\R^{N}_+)=\mu^+:= \frac{N^2}{4}~\!.
\end{equation}

If $N=2$ and if $\C_{\Sigma_\theta}\subset\R^2$ is a cone of amplitude
$\theta\in(0,2\pi]$
then $\lambda_{1}(\Sigma_\theta)$ coincide with  the  Dirichlet  eigenvalue
 on the interval $(0,\theta)$. Hence we get the conclusion, which was first
 pointed out in  \cite{CaMuUMI}:
\begin{equation}
\label{eq:L2-cone} \mu_0(\C_{\Sigma_\theta})=\frac{\pi^2}{\theta^2}\ge
\frac{1}{4}~\!.
\end{equation}

Let $\Sigma$ be a domain in $\S^{N-1}$. If $N\ge 3$ the space
$\mathcal D^{1,2}(\C_\Sigma)$ is defined in a standard way as a
close subspace of $\mathcal D^{1,2}(\R^{N-1})$. Notice that in case
$\Sigma=\S^{N-1}$ it turns out that
$$
\mathcal D^{1,2}(\C_{\S^{N-1}})=\mathcal D^{1,2}(\R^N\setminus\{0\})=
\mathcal D^{1,2}(\R^N)
$$
by a known density result.

If $N=2$ and if $\Sigma$ is properly contained in $\S^1$, then
$\mu_0(\C_\Sigma)>0$ by (\ref{eq:L2-cone}). In this case we can introduce
the space $\mathcal D^{1,2}(\C_\Sigma)$ by completing
$C^\infty_c(\C_\Sigma)$ with respect to the Hilbertian norm
$\left(\int_{\C_\Sigma}|\nabla u|^2~dx\right)^{1/2}$.

The next result is an immediate consequence of the fact that
the Dirichlet eigenvalue problem of $-\Delta$ in the strip $\mathcal
Z_\Sigma$ is never achieved. The same conclusion was already noticed in
\cite{CaMuUMI} in case $N=2$ and in \cite{PT}.

\begin{Proposition}
\label{P:not_achieved_cone}
Let $\Sigma$ be a domain in $\S^{N-1}$. Then $\mu_0(\C_\Sigma)$ is not
achieved in $\mathcal D^{1,2}(\C_\Sigma)$.
\end{Proposition}

\section{Existence}
\label{S:existence} In this Section we show that the condition
$\mu_\lambda(\Omega)< \mu^+=N^2/4$ is sufficient to guarantee the
existence of a minimizer for $\mu_\lambda(\Omega)$. We notice that
throughout this section,  the regularity of $\O$ can be relaxed to
Lipschitz domains which are of class $C^2$ at $0$. We start with a
preliminary result.
\begin{Lemma}\label{lem:Jl1} Let   $\O$ be a smooth domain with  $0\in\partial\O$. Then
$$ \displaystyle{\sup_{\l\in\R}\mu_{\l}(\O) = {\mu^+}}.$$
\end{Lemma}

\proof The proof will be carried out in two steps.

\noindent
\textbf{Step 1.} We claim that $\sup_{\l\in\R}\mu_{\l}(\O) \geq {\mu^+}$.

\noindent
We denote by $\nu$ the interior normal of $\de\O$ at 0. For $\d>0$,
we consider the cone
$$
 \C^{\d}_-:=\left\{
x\in\R^{N-1}~|~x\cdot\nu>-\delta|x|~\right\}.
$$
\noindent Now  fix   $\eps>0$. If $\delta$ is small enough then
$\mu_0(\C^\delta_-)\ge \mu^+-\eps$. Since $\O$ is smooth at $0$ then
there exists a small radius $r>0$ (depending on $\delta$) such that
$\O\cap B_{r_\d}(0)\subset \C^\d_-$.

Next, let  $\psi\in C^\infty(B_{r}(0))$ be a cut-off function, satisfying
$$
0\leq\psi\leq 1~,\quad\psi\equiv 0~\textrm{in}~\R^N\setminus
B_{\frac{r}{2}}(0)~,\quad \psi\equiv 1~\textrm{in}~
B_{\frac{r}{4}}(0)~\!.
$$
We write any $u\in H^1_0(\O)$ as $u=\psi u+(1-\psi)u$, to get
\be\label{eq:uepupumpu} \int_{\O}|x|^{-2}|u|^2~dx\leq
\int_{\O}|x|^{-2}|\psi u|^2~dx+c \int_{\O}|u|^2~dx~\!, \ee
 where the
constant $c$ do not depend on $u$. Since $\psi u \in
\mathcal D^{1,2}(\C^{\d}_-)$ then \be\label{eq:muCduepupumpu}
(\mu^+-\eps)\int_{\O}|x|^{-2}|\psi u|^2~dx\leq {\mu_0( \C^\d_-)}
\int_{\O}|x|^{-2}|\psi u|^2~dx\leq \int_{\O}|\nabla(\psi u)|^2~dx \ee by
our choice of the cone $\C^\d_-$. In addition, we have
$$
\int_{\O}|\nabla(\psi u)|^2~dx\leq \int_{\O}|\nabla u|^2~dx+
\frac{1}{2}\int_{\O}\nabla(\psi^2)\cdot\nabla(u^2)~dx+c \int_{\O} |u|^2~dx~\!.
$$
using integration by parts we get
$$
\int_{\O}|\nabla(\psi u)|^2~dx\leq \int_{\O}|\nabla
u|^2~dx-\frac{1}{2}\int_{\O}\Delta(\psi^2)|u|^2~dx + c\int_{\O} |u|^2~dx.
$$
Comparing with (\ref{eq:uepupumpu}) and (\ref{eq:muCduepupumpu}) we
infer that there exits a positive constant $c$ depending only on
$\delta$ such that
 \be\label{eq:mu0Pos}
  (\mu^+-\eps)\,
{\displaystyle\int_{\O}|x|^{-2}| u|^2~dx}\le {\displaystyle
\int_{\O}|\nabla u|^2~dx+c \displaystyle
\int_{\O}|u|^2~dx}~\!\quad\forall u\in H^1_0(\Omega).
\ee
Hence we get $(\mu^+-\eps)\le \mu_{-c}(\O)$. Consequently
$(\mu^+-\eps)\le\sup_\lambda \mu_\lambda(\O)$, and the conclusion
follows by letting $\eps\to 0$.\\

\noindent
\medskip
\textbf{Step 2:} We claim that $\sup_\lambda \mu_\lambda(\O)\leq \mu^+$.

\noindent
For
$\d>0$ we consider the cone
$$
 \C^{\d}_+:=\left\{
x\in\R^{N-1}~|~x\cdot\nu>\delta|x|~\right\}.
$$
As in the first step, for any   $\delta>0$  there exists $r_\d>0$
such that
 $\C^\d_+\cap B_{r}(0)\subset
\O $  for all $r\in(0,r_\d)$. Clearly by scale invariance, $
\mu_0(\C^\d_+\cap B_{r}(0))= \mu_0(\C^\delta_+)$. For $\e>0$, we
let $\phi\in H^1_0( \C^\d_+\cap B_{r}(0))$ such that
$$
\frac{\displaystyle\int_{\C^\d_+\cap B_{r}(0)}|\nabla \phi|^2~dx}
{\displaystyle\int_{\C^\d_+\cap B_{r}(0)}|x|^{-2}|\phi|^2~dx}
\leq \mu_0(\C^\delta_+)+\e.
$$
From this we deduce that
\begin{eqnarray*}
\mu_{\l}(\O)&\leq& \frac{\displaystyle\int_{\C^\d_+\cap
B_{r}(0)}|\nabla \phi|^2~dx-\l \int_{\C^\d_+\cap
B_{r}(0)}|\phi|^{2}~dx} {\displaystyle\int_{\C^\d_+\cap
B_{r_\d}(0)}|x|^{-2}|\phi|^2~dx}\\
&\leq&  \mu_0(\C^\delta_+) +\e+|\l| ~ \frac{\displaystyle
\int_{\C^\d_+\cap B_{r}(0)}|\phi|^2~dx}
{\displaystyle\int_{\C^\d_+\cap B_{r}(0)}|x|^{-2}|\phi|^2~dx}.
\end{eqnarray*}
Since $\displaystyle\int_{\C^\d_+\cap B_{r}(0)}|x|^{-2}|\phi|^2~dx \geq r^{-2}
\int_{\C^\d_+\cap B_{r}(0)}|\phi|^2~dx $, we get
$$
\mu_{\l}(\O)\leq   \mu_0(\C^\delta_+)+\e+ r^2|\l|.
$$
The conclusion  follows immediately,  since $
\mu_0(\C^\delta_+)\to\mu^+$ when $\d\to0$.
 \QED

Notice that if $\O$ is bounded then by \eqref{eq:mu0Pos} and
Poincar\'e inequality
\begin{equation}
\label{eq:mu0positive} \mu_0(\O)>0~\!.
\end{equation}
For $N=2$ this  was shown in \cite{CaMuUMI} and  for   more general
domains. We are in position to prove the main result of this
section.
\begin{Theorem} \label{th:extJl0}
 Let $\lambda\in \R$ and let $\O$ be a smooth bounded domain  of $\R^N$ with $0\in\partial\Omega$.
If  $\mu_\l(\O)<\mu^+$ then $\mu_\l(\O)$  is attained.
\end{Theorem}

\proof
Let $u_n\in H^1_0(\O)$ be a minimizing sequence for $\mu_\lambda(\O)$.
We can normalize it to have
\begin{gather}
\label{eq:normun}
\int_{\O}|\nabla u_n|^2=1,\\
\label{eq:nutoJpo} 1-\l \int_{\O}| u_n|^2=\mu_\l(\O)\int_\O|x|^{-2}|u_n|^2+o(1)~\!.
\end{gather}
We can assume that $u_n\rightharpoonup u$ weakly in $H^1_0(\O)$,
$|x|^{-1}u_n \rightharpoonup |x|^{-1}u$ weakly in $L^2(\O)$,
and
$u_n  \to u $ in $L^2(\O)$, by (\ref{eq:mu0positive}) and by Rellich Theorem.
Putting $\theta_n:= u_n-u$, from  \eqref{eq:normun} and \eqref{eq:nutoJpo} we get
\begin{gather}
\nonumber
\int_{\O}|\nabla \theta_n|^2+\int_\O|\nabla u|^2=1+o(1),\\
\label{eq:JunearJ}
1-\l\int_{\O}|u|^2=\mu_\l(\O)\left(\int_\O|x|^{-2}|\theta_n|^2+\int_\O|x|^{-2}|u|^2\right)+o(1)~\!.
\end{gather}
By Lemma~\ref{lem:Jl1},
 for any fixed positive $\d<\mu^+-\mu_\lambda(\O)$, there exists $\l_\d\in\R$ such that $\mu_{\lambda_\delta}(\O)\ge \mu^+-\delta$.
 Hence
 $$
\int_\O|\nabla \theta_n|^2+o(1)\ge
(\mu^+-\delta)\int_\O|x|^{-2}|\theta_n|^2~\!,
 $$
 as $\theta_n\to 0$ in $L^2(\O)$.  Testing $\mu_\lambda(\O)$ with $u$ we get
 \begin{eqnarray*}
 \mu_\lambda(\O)\int_\O|x|^{-2}|u|^2&\le&
 \int_\O|\nabla u|^2-\l \int_{\O}|u|^2\le 1-\int_{\O}|\nabla \theta_n|^2-\l \int_{\O}|u|^2+o(1)\\
 &\le&1-(\mu^+-\delta)\int_\O|x|^{-2}|\theta_n|^2-\l \int_{\O}|u|^2+o(1)\\
&\le&(\mu_\lambda(\O)-\mu^++\delta)\int_\O|x|^{-2}|\theta_n|^2+
\mu_\lambda(\O)\int_\O|x|^{-2}|u|^2+o(1)
\end{eqnarray*}
by (\ref{eq:JunearJ}). Therefore $\int_\O|x|^{-2}|\theta_n|^2\to 0$, since
$\mu_\lambda(\O)-\mu^++\delta<0$. In particular,
$$
\mu_\lambda(\O)\int_\O|x|^{-2}|u|^2= \int_\O|\nabla u|^2-\l
\int_{\O}|u|^2
$$
and  $u\neq 0$ by (\ref{eq:JunearJ}). Thus
$u$ achieves $\mu_\lambda(\O)$.
\QED

We conclude this section with a corollary of Theorem
\ref{th:extJl0}.

\begin{Corollary}
 Let $\O$ be a smooth bounded domain  of $\R^N$ with $0\in\partial\Omega$.
 Then
 $$
 \frac{(N-2)^2}{4}<\mu_0(\O)\le \frac{N^2}{4}~\!.
 $$
\end{Corollary}

\proof
It has been already proved in Lemma \ref{lem:Jl1} that
$\mu_\lambda(\O)\le \frac{N^2}{4}$. If the strict inequality holds, then
there exists $u\in H^1_0(\O)$ that achieves $\mu_0(\O)$, by Theorem
\ref{th:extJl0}. But then  $\frac{(N-2)^2}{4}<\mu_0(\O)$, otherwise
a null extension of $u$ outside $\O$ would achieve the Hardy constant on $\R^N$.
\QED

\begin{Remark}
Following \cite{CaMuUMI}, for non smooth domains $\O$ we can
introduce the "limiting" Hardy constant
$$
\hat \mu_0(\Omega)=\sup_{r>0}\mu_0(\Omega\cap B_{r})~\!.
$$
Using similar arguments it can be proved that $\sup_{\l}\mu_{\l}(\O)
= \hat \mu_0(\Omega)$, and that $\mu_\l(\O)$ is achieved provided
$\mu_\l(\O)<\hat \mu_0 (\O)$.
\end{Remark}

\section{Nonexistence}
\label{S:ne}

The main result in this section is stated in the following theorem.

 \begin{Theorem}
\label{T:ne} Let $\O$ be a domain in $\R^N$, $N\ge 2$, and let $\lambda\in\R$.
Assume that
there exist $R>0$ and a Lipschitz domain $\Sigma\subset \S^{N-1}$
such that $B_R\cap \C_\Sigma\subset\O$. If $u\in \widehat H^1(\O)$
solves
\begin{equation}
\label{eq:pbne}
\begin{cases}
-\Delta
u\ge\displaystyle{\left(\frac{(N-2)^2}{4}+\lambda_1(\Sigma)\right)}
|x|^{-2} u+\lambda u
&\textit{ in $\mathcal D'(\Omega\setminus\{0\})$}\\
u\ge 0~\!,
\end{cases}
\end{equation}
then $u\equiv 0$ in $\O$.
\end{Theorem}

Before proving Theorem \ref{T:ne} we point out some of its
consequences.

\begin{Corollary}
\label{C:ne} Let $\Omega$ be a smooth bounded domain
containing a half-ball and such that $0\in\partial\O$. If
$\mu_\lambda(\O)=\mu^+$ then $\mu_\lambda(\O)$ is not achieved.
\end{Corollary}

\proof
Assume that $u$ achieves $\mu_\l(\O)=\mu^+$. Then $u$ is a weak solution to
\begin{equation}
\label{eq:corollary}
-\Delta u = \mu^+ |x|^{-2} u+\lambda u~\!.
\end{equation}
Test (\ref{eq:corollary}) with the negative and the positive part of
$u$ to conclude that $u$ has constant sign. Now by the maximum
principle  $u>0$ in $\O$, contradicting Theorem \ref{T:ne}, since
$\O\supset B_R\cap \C_{\S_+^{N-1}}$ and $\lambda_1(\S_+^{N-1})=N-1$.
\QED

We also point out the following consequence to Theorem \ref{T:ne}, that holds
for smooth
domains $\O$ with $0\in\partial\O$.

\begin{Theorem}
\label{T:ne2} Let $\O$ be a smooth domain in $\R^N$, $N\ge 2$ with
$0\in\partial\O$
 and let $\lambda\in\R$. If $u\in \widehat H^1(\O)$
solves
$$
\begin{cases}
-\Delta u\ge\mu |x|^{-2} u+\lambda u
&\textit{ in $\mathcal D'(\Omega)$}\\
u\ge 0~\!
\end{cases}
$$
for some $\mu>\mu^+$, then $u\equiv 0$ in $\O$.
\end{Theorem}

\proof We start by noticing that there exists a geodesic ball
$\Sigma\subset \S^{N-1}$ contained in a hemisphere, and such that
$\lambda_1(\Sigma)\le N-1+\mu-\mu^+$. Since $0\in\de\O$ and since
$\partial\O$ is smooth then, up to a rotation, we can find a small
radius $r>0$ such that $B_{r}\cap \C_{\Sigma}\subset \O$. The
conclusion follows from Theorem \ref{T:ne}, as $\displaystyle{\mu
\ge (N-2)^2/4+\lambda_1(\Sigma)}$. \QED

\begin{Remark}
Theorem \ref{T:ne} applies also when the origin lies in the interior of the domain. More precisely, let $\O$ be any domain in $\R^N$, with $N\ge 2$
and $0\in\O$. If $u\in \widehat H^1_{\rm loc}(\O)$ is a nonnegative solution
to
$$
-\Delta u\ge \displaystyle{\frac{(N-2)^2}{4}}~|x|^{-2} u+\lambda u
\quad\textit{ in $\mathcal D'(\Omega\setminus\{0\})$}
$$
for some $\lambda\in\R$, then $u\equiv 0$ in $\O$.
\end{Remark}

\bigskip

In order to prove Theorem \ref{T:ne} we need few preliminary results
about maps of two variables. Recall that $\disc_R\subset \R^2$ is
the open disk of radius $R$ centered at $0$.

\begin{Lemma}
\label{L:psi} Let $\psi\in\widehat H^1(\disc_R)$ and $f\in
L^1_{\textrm{loc}}(\disc_R)$ for some $R>0$. If $\psi$ solves
\begin{equation}
\label{eq:psi1} -\Delta \psi\ge f \quad \textit{ in
$\mathcal D'(\disc_R\setminus\{0\})$}
\end{equation}
then $-\Delta \psi\ge f$ in $\mathcal D'(\disc_R)$.
\end{Lemma}

\proof We start by noticing that from
$$
\infty>\int_{\disc_R}|z|^{-2}|\psi|^2=\int_0^R\frac{1}{r}\left(r^{-1}\int_{\partial
B_r}|\psi|^2\right)
$$
it follows that there exists a sequence $r_h\to 0$, $r_h\in(0,R)$
such that
\begin{equation}
\label{eq:rh} r_h^{-1}\int_{\partial B_{r_h}}|\psi|^2\to 0~,\quad
r_h^{-2}\int_{\partial B_{r_h^2}}|\psi|^2\to 0
\end{equation}
as $h\to \infty$. Next we introduce the following cut-off functions:
$$
\eta_h(z)=\begin{cases}
0&\textrm{if $|z|\le r_h^2$}\\
\displaystyle{\frac{\log |z|/r_h^2}{|\log r_h|}}&\textrm{if $r_h^2<|z|<r_h$}\\
1&\textrm{if $r_h\le |z|\le R$.}
\end{cases}
$$
Let $\f\in C^\infty_c(\disc_R)$ be any nonnegative function. We test
(\ref{eq:psi1}) with $\eta_h\f$  to get
$$
\int\nabla \psi\cdot\nabla(\eta_h\f)\ge\int f~\!\eta_h\f~.
$$
Since $\psi\in H^1(\disc_R)$ and since $\eta_h\weak 1$ weakly$^*$ in
$L^\infty$, it is easy to check that
$$
\int f~\!\eta_h\f=\int f\f+o(1)~,\quad
\int\eta_h\nabla\psi\cdot\nabla\f=\int\nabla\psi\cdot\nabla\f+o(1)
$$
as $h\to \infty$. Therefore
\begin{equation}
\label{eq:test}
\int\nabla \psi\cdot\nabla\f+\int\f\nabla\psi\cdot\nabla\eta_h\ge\int f~\!\f+o(1)~\!.%
\end{equation}
To pass to the limit in the left-hand side we notice that
$\nabla\eta_h$ vanishes outside the annulus
$A_h:=\{r_h^2<|z|<r_h\}$, and that $\eta_h$ is harmonic on $A_h$.
Thus
\begin{eqnarray*}
\int\f\nabla\psi\cdot\nabla\eta_h&=&\int_{A_h}\nabla(\psi\f)\cdot\nabla\eta_h-
\int_{A_h}\psi\nabla\f\cdot\nabla\psi\\
&=&\mathcal R_h - \int_{A_h}\psi\nabla\f\cdot\nabla\eta_h
\end{eqnarray*}
where
$$
\mathcal R_h:=-r_h^{-2}\int_{\partial B_{r_h^2}}(\nabla \eta_h\cdot
z)\psi\f+ r_h^{-1}\int_{\partial B_{r_h}}(\nabla \eta_h\cdot z
)\psi\f~\!.
$$
Now
$$
|\mathcal R_h|\le c~\!(r_h|\log r_h|)^{-1}\int_{\partial
B_{r_h}}|\psi|+ c~\!(r^2_h|\log r_h|)^{-1}\int_{\partial
B_{r^2_h}}|\psi|
$$
where $c>0$ is a constant that does not depend on $h$, and
$$
(r_h|\log r_h|)^{-1}\int_{\partial B_{r_h}}|\psi|\le c~\!|\log
r_h|^{-1}\left(r_h^{-1}\int _{\partial
B_{r_h}}|\psi|^2\right)^{1/2}=o(1)
$$
by H\"older inequality and by (\ref{eq:rh}). In the same way, also
$$
(r^2_h|\log r_h|)^{-1}\int_{\partial B_{r^2_h}}|\psi|\le c~\!|\log
r_h|^{-1}\left(r^{-2}_h\int _{\partial
B_{r^2_h}}|\psi|^2\right)^{1/2}=o(1)~\!,
$$
and hence $\mathcal R_h=o(1)$.
 Moreover from
$\psi\in L^2(\disc_R;|z  |^{-2}~dz  )$ it follows that
$$
\left|\int_{A_h}\psi\nabla\f\cdot\nabla\eta_h\right|\le |\log
r_h|^{-1}\int|z  |^{-1}\psi|\nabla\f|= o(1)~.
$$
In conclusion, we have proved that
$\displaystyle\int\f\nabla\psi\cdot\nabla\eta_h= o(1)$ and therefore
(\ref{eq:test}) gives
$$
\int\nabla \psi\cdot\nabla\f\ge \int f~\!\f~\!.
$$
Since $\f$ was an arbitrary nonnegative function in
$C^\infty_c(\disc_R)$, this proves that $-\Delta\psi\ge f$ in
the distributional sense on $\disc_R$, as desired. \QED

\noindent The same proof gives a similar result for subsolutions.

\begin{Lemma}
\label{L:psi2} Let $\f\in\widehat H^1(\disc_R)$ and $f\in
L^1_{\textrm{loc}}(\disc_R)$ for some $R>0$. If $\f$ solves
$$
\Delta \f\ge f \quad \textit{ in $\mathcal
D'(\disc_R\setminus\{0\})$}
$$
then $\Delta \f\ge f$ in $\mathcal D'(\disc_R)$.
\end{Lemma}

The next result is crucial in our proof. We state it is a more general form than needed, as
it could have an independent interest. Notice that we do not need any a priori knowledge of the
sign of $\psi$ in the interior of its domain.

\begin{Lemma}
\label{P:psi} For any $\lambda\in \R$ there exists $R_\lambda>0$ such that for any $R\in(0,R_\lambda)$, $\eps>0$,
problem
\begin{equation}
\label{eq:pbneP}
\begin{cases}
-\Delta \psi\ge \lambda \psi
&\textit{ in $\mathcal D'(\disc_R\setminus\{0\})$}\\
\psi\ge \eps&\textrm{on $\partial \disc_R$.}
\end{cases}
\end{equation}
has no solution $\psi\in \widehat H^1(\disc_R)$.
\end{Lemma}

\proof We fix $R_\lambda<1/3$ small enough, in such a way that
\begin{equation}
\label{eq:Rlambda2} \lambda<\lambda_1(\mathbb D_{R_\lambda})~\quad
\textrm{if $\lambda\ge 0$}~\!,
\end{equation}
\begin{equation}
\label{eq:Rlambda}
|\lambda||z|^2\left|\log|z|\right|^2\le\frac{3}{4}~\quad
\textrm{for any $z\in\disc_{R_\lambda}$~,~if
$\lambda< 0$}~\!,
\end{equation}
We claim that the conclusion in Lemma \ref{P:psi} holds with this
choice of $R_\lambda$. We argue by contradiction. Let $R<R_\lambda$
and $\eps>0$, $\psi\in \widehat H^1(\mathbb D_R)$ as in
(\ref{eq:pbneP}).

For any $\delta\in(1/2,1)$ we introduce the following radially
symmetric function on $\mathbb D_R$:
$$
\f_\delta(z)=\left|\log |z|\right|^{-\delta}~.
$$
By direct computation one can easily check that
$\f_\delta\in\widehat H^1(\mathbb D_R)$, and in particular
\begin{equation}
\label{eq:L2phi} (2\delta-1)\int_{\mathbb
D_R}|z|^{-2}|\f_\delta|^2=2\pi+o(1)\quad\textrm{as $\delta\to
\frac{1}{2}$.}
\end{equation}
Since $\delta>1/2$ then $\f_\delta$ is a smooth solution to
\begin{equation}
\label{eq:Delta-phi}
\Delta\f_\delta\ge\frac{3}{4}~\!|z|^{-2}\left|\log|z|\right|^{-2+\delta}=
\frac{3}{4}~\!|z|^{-2}\left|\log|z|\right|^{-2}\f_\delta
\end{equation}
in $\mathbb D_R\setminus\{0\}$. By Lemma \ref{L:psi2} we infer that
$\f_\delta$  solves (\ref{eq:Delta-phi}) in the dual of $\widehat
H^1(\mathbb D_R)$. Next we put
$$
v:=\eps\f_\delta-\psi\in \widehat H^1(\mathbb D_R)~\!,
$$
and we notice that $v\le 0$ on $\partial \mathbb D_R$, as $R<1/3$.
Notice also that
\begin{eqnarray*}
\Delta v&\ge& \frac{3}{4}~\!|z|^{-2}\left|\log|z|\right|^{-2}(\eps\f_\delta)+ \lambda\psi\\
&=&\left[\frac{3}{4}~\!|z|^{-2}\left|\log|z|\right|^{-2}+\lambda\right](\eps\f_\delta)-\lambda
v
\end{eqnarray*}
on  the dual of $\widehat H^1(\mathbb D_R)$, by (\ref{eq:Rlambda}). We use as
test function $v^+:=\max\{v,0\}\in H^1_0(\mathbb D_R)\cap \widehat
H^1(\mathbb D_R)$ to get
$$
-\int_{\mathbb D_R}|\nabla v^+|^2\ge \int_{\mathbb D_R}
\left[\frac{3}{4}~\!|z|^{-2}\left|\log|z|\right|^{-2}+\lambda\right](\eps\f_\delta)v^+-
\lambda\int_{\mathbb D_R}|v^+|^2~\!.
$$
If $\lambda\ge 0$ we infer that
$$
\int_{\mathbb D_R}|\nabla v^+|^2\le \lambda\int_{\mathbb D_R}|v^+|^2
$$
and hence $v^+\equiv 0$ on $\mathbb D_R$ by (\ref{eq:Rlambda2}). If
$\lambda<0$ we get
$$
0\ge-\int_{\mathbb D_R}|\nabla v^+|^2\ge|\lambda|\int_{\mathbb
D_R}|v^+|^2~\!,
$$
hence again $v^+=0$ on $\mathbb D_R$, by (\ref{eq:Rlambda}). Thus $\psi\ge \eps \f_\delta$ on
$\mathbb D_R$ and therefore
$$
\infty>\int_{\mathbb D_R}|z|^{-2}|\psi|^2\ge \eps \int_{\mathbb
D_R}|z|^{-2}|\f_\delta|^2~,
$$
which contradicts (\ref{eq:L2phi}). \QED

\bigskip

\noindent {\bf Proof of Theorem \ref{T:ne}.} Without loss of
generality, we may assume that $\lambda<0$. Let $\Phi>0$ be the
first eigenfunction of $-\Delta_\sigma$ on $\Sigma$. Thus $\Phi$
solves \be\label{eq:Phi-Sig}
\begin{cases}
-\Delta_\sigma\Phi=\lambda_1(\Sigma)\Phi&\textrm{in $\Sigma$}\\
\Phi=0~,\quad \displaystyle\frac{\partial\Phi}{\partial\eta}\le 0&
\textrm{on $\partial\Sigma$,}
\end{cases}
\ee where $\eta\in T_\sigma(\S^{N-1})$ is the exterior normal to
$\Sigma$ at $\sigma\in\partial \Sigma$.

 By density and the trace theorem, we can define the radially
symmetric map $\psi$ in $\disc_R\setminus\{0\}$ as
 \begin{equation}
 \label{eq:defpsi}
 \psi(z)=|z|^{\frac{N-2}{2}}\int_{\Sigma} u(|z|\sigma)\Phi(\sigma)~d\sigma~\!
 =|z|^{\frac{N-2}{2}}\int_{|z|\Sigma} u(\s')\Phi_{|z|}(\s')~d\s'\!,
 \end{equation}
where $\Phi_r(\s')=\Phi(\frac{\s'}{r})$ for all $\s'\in r\Sigma$.
Since in polar coordinates $(r,\sigma)\in(0,\infty)\times\S^{N-1}$
it holds that
 $$
 u_{rr}=-(N-1)r^{-1}u_r-r^{-2}\Delta_\sigma u~,
 $$
 direct computations based on \eqref{eq:pbne} lead to
 $$
 -\Delta\psi\ge \lambda\psi ~~\textrm{in $\mathcal
D'(\disc_R\setminus\{0\})$}.
$$
 We claim that $\psi\in \widehat H^1(\disc_R)$. Indeed, for $r=|z|$,
 $$
 |\psi'|\le c r^{\frac{N-2}{2}-1}\int_{\Sigma}|u(r\sigma)| +
 c r^{\frac{N-2}{2}}\int_{\Sigma}|\nabla u(r\sigma)| ~,
 $$
 and, by H\"older inequality,
$$
\int_{\disc_R}\left(r^{\frac{N-2}{2}-1}\int_{\Sigma}|u(r\sigma)|\right)^2=
c\int_0^R\int_{\Sigma}r^{N-3}u^2\le c\int_\O|x|^{-2}u^2<\infty~,
$$
$$
\int_{\disc_R}\left( r^{\frac{N-2}{2}}\int_{\Sigma}|\nabla
u(r\sigma)|\right)^2 \le c\int_0^Rr^{N-1}\int_\Sigma|\nabla u|^2\le
c \int_\O|\nabla u|^2<\infty.
$$
Finally, $\psi\in L^2(R^2_R;|z|^{-2}dz)$ as
$$
\int_{\disc_R}|z|^{-2}|\psi|^2=2\pi\int_0^Rr^{-1}|\psi|^2 \le
c\int_0^R r^{N-3}\int_\Sigma|u|^2=c\int_{\O}|x|^{-2}|u|^2<\infty~\!.
$$
Thus Lemma \ref{P:psi} applies and
 since $\psi$ is radially symmetric we get $\psi\equiv 0$ in a
neighborhood of $0$. Hence $u\equiv 0$ in $B_r\cap\C_\Sigma$, for $r>0$ small enough.
To conclude the proof in case $\Omega$ strictly contains $B_r\cap\C_\Sigma$, take
any domain $\O'$ compactly contained in $\O\setminus\{0\}$ and such that
$\O'$ intersects $B_r\cap\C_\Sigma$. Via a convolution procedure, approximate
$u$ in $H^1(\O')$ by a sequence of smooth maps $u_\eps$ that solve
$$
-\Delta u_\eps+ |\lambda|u_\eps\ge 0\quad\textrm{in $\Omega'$}~\!.
$$
Since $u_\eps\ge 0$ and $u_\eps\equiv 0$ on $\O'\cap
B_r\cap\C_\Sigma$, then $u_\eps\equiv 0$ on $\O'$ by the  maximum
principle. Thus also $u\equiv 0$ in $\O'$, and the conclusion
follows. \QED

\section{Remainder terms}
\label{S:reminder}

We prove here some inequalities that will be used in the next section to estimate the infimum
$\lambda^*$ defined in (\ref{eq:lambda*}).

Brezis and V\'{a}zquez proved in \cite{BV} the following improved Hardy inequality:
\begin{equation}
\label{eq:BV}
\int_\Omega|\nabla u|^2-\frac{(N-2)^2}{4}\int_\Omega|x|^{-2}|u|^2\ge
\omega_N\frac{\lambda(\disc)}{|\Omega|}\int_\Omega|u|^2~\!,
\end{equation}
that holds for any $u\in C^\infty_c(\Omega)$. Here $\Omega\subset\R^N$ is any bounded domain,
$\lambda(\disc)$ is the first Dirichlet eigenvalue of the unit ball
$\disc$ in $\R^2$, and $\omega_N$, $|\Omega|$ denote the measure of the
unit ball in $\R^N$ and of $\Omega$, respectively.
If $0\in\Omega$ then $(N-2)^2/4$ is the Hardy constant $\mu_0(\Omega)$ relative to the
domain $\Omega$,  by the invariance of the ratio
$$
\frac{\displaystyle\int_{\O}|\nabla u|^2~dx}
{\displaystyle\int_{\O}|x|^{-2}|u|^2~dx}~
$$
with respect to dilations in $\R^N$.

We show that a Brezis-V\'{a}zquez type inequality holds  in case the
singularity is placed at the boundary of the domain. We start with
conic domains
$$
\C_{R,\Sigma}=\{ t\sigma~|~t\in(0,R),\, \sigma\in \Sigma~\},
$$
where $\Sigma\subset \S^{N-1}$ and $R>0$.

\begin{Proposition}
\label{lem:reminder_cone}
Let $\Sigma$ be a domain in  $\S^{N-1}$. Then
\begin{equation}
\label{eq:reminder_cone} \int_{\C_{R,\Sigma}}
|\nabla u|^2-\mu_0\left(\C_{\Sigma}\right)\,\int_{\C_{R,\Sigma}}|x|^{-2}|u|^2\geq
\frac{\lambda_1(\disc)}{R^2}\,\int_{\C_{R,\Sigma}}|u|^2,\qquad\forall u\in
C^{\infty}_c(\C_{1,\Sigma}).
\end{equation}
\end{Proposition}

\proof
By homogeneity, it suffices to prove the proposition for $R=1$.
Fix $u\in C^{\infty}_c(\C_{1,\Sigma})$ and compute in  polar coordinates $t=|x|$, $\sigma=x/|x|$:
$$
\int_{\C_{1,\Sigma}}|\nabla u|^2= \int_0^1\int_\Sigma\left|\frac{\de u}{\de t}\right|^2t^{N-1}dtd\s+
 \int_0^1\int_\Sigma\left|\nabla_\s u\right|^2t^{N-3}dtd\s~\!,
$$
$$
\int_{\C_{1,\Sigma}}|x|^{-2}|u|^2=\int_0^1\int_\Sigma |u|^2
t^{N-3}dtd\s.
$$
 Since for every $t\in(0,1)$ it holds that
$$
\int_\Sigma\left|\nabla_\s u\right|^2 t^{N-3}d\s \geq
\lambda_1(\Sigma)\,\int_\Sigma |u|^2 t^{N-3}d\s,
$$
then by Proposition \ref{P:mu} we only have to show that
\be\label{eq:imphc} \int_0^1\left|\frac{\de u}{\de
t}\right|^2t^{N-1}dt-\frac{(N-2)^2}{4}\, \int_0^1 |u|^2
t^{N-3}dt\geq \lambda_1(\disc)\,\int_0^1 |u|^2 t^{N-1}dt \ee
for any fixed $\s\in \Sigma$.
For that, we put
$w(t)=t^{\frac{N-2}{2}}u(t\sigma)$, and we compute
\begin{eqnarray*}
\int_0^1\left|\frac{\de u}{\de
t}\right|^2t^{N-1}dt&-&\mu_0\left(\R^N\right)\, \int_0^1 |u|^2
t^{N-3}dt\\
&=&
\int_0^1\left|\frac{\de w}{\de t}\right|^2t dt+\left(2-N\right)\, \int_0^1 \frac{\de w}{\de t}wdt\\
&=&
\int_0^1\left|\frac{\de w}{\de t}\right|^2t dt+
\frac{\left(2-N\right)}{2}\,\int_0^1\frac{\de w^2}{\de t} dt
=\int_0^1\left|\frac{\de w}{\de t}\right|^2t dt\\
&\geq&\lambda_1(\disc)\,\int_0^1w^2t  dt=\lambda_1(\disc)\,\int_0^1 |u|^2 t^{N-1}dt.
\end{eqnarray*}
This gives \eqref{eq:imphc} and the proposition  is proved.
\QED

The main result in this section is contained in the next theorem.

\begin{Theorem}\label{th:bd-ls}
Let $\O$ be a bounded domain of $\R^N$ with $0\in\de\Omega$. If $\O$
is contained in a half-space then
$$
\int_{\O}|\nabla u|^2-\mu^+\int_{\O} |x|^{-2}|u|^2\geq
\frac{\lambda_1(\disc)}{|\textrm{diam}(\O)|^2  }\int_{\O}
|u|^2\qquad \forall u\in H^1_0(\O).
$$
\end{Theorem}

\proof Let $R>0$ be the diameter of $\Omega$. Then $\Omega\subset
B^+_R$, where $B^+_R$ is a half  ball of radius $R$  centered at the
origin. Take $\Sigma$ to be a half sphere in $\S^{N-1}$  in
Proposition \ref{lem:reminder_cone}, so that $\C_\Sigma$ is an
half-space. Recalling  (\ref{eq:half-space}),  we conclude that
$$
\int_{B^+_R}|\nabla u|^2-\mu^+\int_{B^+_R} |x|^{-2}|u|^2\geq
\frac{\lambda_1(\disc)}{R^2}\int_{B^+_R}| u|^2
$$
for any $R>0$, $u\in C^{\infty}_c(\O)$ and   the theorem readily
follows. \QED

\begin{Remark}
\label{P:reminder2}
Let $\O$ be a bounded domain of $\R^2$ with
$0\in\de\Omega$ and assume that $\Omega$ does not intersect an half-line
emanating from the origin. Then (\ref{eq:L2-cone}) and
Proposition \ref{lem:reminder_cone} imply the following
  improved Hardy inequality:
$$
\int_{\O}|\nabla u|^2-\frac{1}{4}\int_{\O} |x|^{-2}|u|^2\geq
\frac{\lambda_1(\disc)}{|\textrm{diam}(\O)|^2  }\int_{\O} |u|^2\qquad \forall
u\in H^1_0(\O).
$$
\end{Remark}

\begin{Remark}
\label{rem:Lp}
 As pointed out by Brezis-V\'{a}zquez in~\cite{BV}, Extension 4.3, the following
 Hardy-Sobolev inequality holds
$$
\int_{\C_{1,\Sigma}}|\nabla
u|^2-\mu_0\left(\C_{1,\Sigma}\right)\,\int_{\C_{1,\Sigma}}|x|^{-2}|u|^2\geq
c_p\,\left(\int_{\C_{1,\Sigma}}|u|^p\right)^{\frac{2}{p}},\qquad\forall
u\in C^{\infty}_c(\C_{1,\Sigma})
$$
for all $p\in \left(2,\frac{2N}{N-2}\right)$, where $c_p$ is a
positive constant depending on $p$ and $N$.
\end{Remark}

\section{Estimates on $\lambda^*$}
\label{S:estimates}

In this section we provide sufficient conditions to have $\lambda^*>-\infty$
or
$\lambda^*<0$.

\subsection{Estimates from below}

Let $\O$ be
a smooth domain in $\R^N$ with $0\in\partial \O$. We say that
$\O$ is {\em locally  convex} at $0$ if exists a ball $B$
centered at $0$ such that $\Omega\cap B$ is contained in a
half-space. In essence, for domains of class $C^2$ this means
that all the principal
curvatures of $\de\O$ (with respect to the interior normal) at 0 are
strictly positive.

In the case where $\O$ is locally convex at $0\in\de\O$, the
supremum in Lemma~\ref{lem:Jl1} is attained.

 \begin{Proposition}\label{lem:Jlconvx}
If $\O$ is
locally convex at $0$, then there exists $\l^*(\O)\in\R$ such that
$$
\begin{array}{cc}
  \mu_{\l}(\O)=\mu^+, &  \quad\forall
\l\leq\l^*(\O), \\
 \mu_{\l}(\O)<\mu^+,  & \quad\forall
\l>\l^*(\O).
\end{array}
$$
 \end{Proposition}
 \proof
 The locally convexity assumption at $0$ means that there
exists $r>0$ such that $ B_r(0)\cap\O$ is contained in a half space.
We let $\psi \in C^\infty_c(\R^N)$ with $0\leq \psi\leq 1$,
$\psi\equiv 0$ in $\R^N\setminus B_\frac{r}{2}(0)$ and
 $\psi\equiv 1$ in $B_{\frac{r}{4}}(0)$.
 Arguing as in the proof of Lemma~\ref{lem:Jl1}, for every  $u\in H^1_0(\O)$
 we get
\be\label{eq:uepupumpu1} \int_{\O}|x|^{-2}|u|^2~dx\leq
\int_{\O}|x|^{-2}|\psi u|^2~dx+c\int_{\O}|u|^2~dx \ee
for some constant $c=c(r)>0$.
Since $\psi u \in H^1_0(B_r(0)\cap\O)$, from the definition of $\mu^+$
we infer
$$
{\mu^+} \int_{\O}|x|^{-2}|\psi u|^2~dx\leq     \int_{\O}|\nabla(\psi u)|^2~dx.
$$
As in Lemma \ref{lem:Jl1} we get
$$
\int_{\O}|\nabla(\psi u)|^2~dx\leq \int_{\O}|\nabla u|^2~dx+c \int_{\O}~\!
|u|^2~dx~\!.
$$
Comparing with (\ref{eq:uepupumpu1}), we infer that there exits a positive
constant $c$  such that
$$
{\mu^+}\int_{\O}|x|^{-2}|u|^2~dx\leq   \int_{\O}|\nabla u|^2~dx+
c\int_{\O} |u|^2~dx.
$$
This proves that $\mu_{-c}(\O)\geq\mu^+$. Thus $\mu_{-c}(\O)=\mu^+$ by Lemma \ref{lem:Jl1}.
 Finally, noticing that  $\mu_\l(\O)$ is
  decreasing in $\l$, we can set
\be \l^*(\O):=\sup\{{\l\in\R}\quad:\quad \mu_\l(\O)= {\mu^+}\}
 \ee
 so that $\mu_\l(\O)<\mu^+$ for all $\l>\l^*(\O)$.
\QED

Finally, we notice that by Lemma~\ref{lem:Jl1}, if  $\O$ is contained
in a half-space, then $\mu_0(\O)=\mu^+$, and
therefore $\l^*(\O)\geq0$. Thus, From  Theorem~\ref{th:bd-ls} we infer the following result.

\begin{Theorem}
\label{rem:lsneg3}
Let $\O$ be a bounded smooth domain with $0\in\partial\O$. If
$\O$ is contained in a half-space then
$$
\l^*(\O)\geq \frac{\lambda_1(\disc)}{|\textrm{diam}(\O)|^2}~\!.
$$
\end{Theorem}

It would be of interest to know if it is possible to get lower bounds
depending only on the measure of $\Omega$,  as in \cite{BV} and \cite{HO}.

\subsection{Estimates from above}

The  local convexity assumption of $\O$ at $0$ does not necessary
implies that $\l^*(\O)\geq 0$. Indeed the following remark holds.

\begin{Proposition}
\label{P:fat}
For any $\delta>0$ there exists $\rho_\delta>0$ such that if $\O$ is a smooth domain
with $0\in\partial\O$ and
$$
\O\supseteq \{x\in\R^N~|~x\cdot\nu>-\delta|x|~,~\alpha<|x|<\beta~\}
$$
for some $\nu\in\S^{N-1}$, $\beta>\alpha>0$ with $\beta/\alpha>\rho_\delta$, then
$\lambda^*<0$. In particular the Hardy constant $\mu_0(\O)$ is achieved.
\end{Proposition}

\proof
Since the cone
$$
\C_\delta=\{ x\in\R^N~|~x\cdot\nu>-\delta|x|~\}
$$
contains a hemispace, then its Hardy constant is smaller than $\mu^+$. Thus
there exists $u\in C^\infty_c(\C_\delta)$ such that
$$
\frac{\displaystyle\int_{\C_\delta}|\nabla u|^2~dx}
{\displaystyle\int_{\C_\delta}|x|^{-2}|u|^2~dx} <\mu^+~.
$$
Assume that the support of $u$ is contained in an annulus of radii
$b>a>0$. Then the conclusion in Proposition \ref{P:fat} holds, with  $\rho:=b/a$.
\QED

Notice that $\O$  can  be locally strictly convex at
 $0$.

\begin{Remark}
\label{R:GR}
A similar remark holds for the following minimization problem,
which is related to the Caffarelli-Kohn-Nirenberg inequalities:
\begin{equation}
\label{eq:GR}
\inf_
{\genfrac{}{}{0pt}{}{\scriptstyle{u\in H^1_0(\O)}}{\scriptstyle{u\neq 0}}}
~\frac{\displaystyle\int_{\O}|\nabla u|^2~dx}
{\displaystyle\left(\int_{\O}|x|^{-b}|u|^p~dx\right)^{2/p}}~\!,
\end{equation}
where $2<p<2^*$, $b:=N-p(N-2)/2$. In case $0\in\partial\O$, the minimization problem
(\ref{eq:GR}) was studied in
\cite{GK}, \cite{GR1} and \cite{GR2}.
\end{Remark}

\begin{Remark}\label{rem:lsneg2}
We do not know wether the strict local concavity of $\O$ at 0 can implies that
$\mu_0(\O)<\mu^+$. See the paper \cite{GK} by  Ghoussoub and Kang  for the minimization
problem (\ref{eq:GR}).
\end{Remark}

\label{References}

\end{document}